\documentclass{article}
\title{On Isomorphisms between Certain Yetter-Drinfel'd Hopf Algebras}
\author{Yevgenia Kashina, Yorck Sommerh\"auser}
\date{}

\usepackage{amsmath}
\usepackage{theorem}
\usepackage{amssymb}
\usepackage{amscd}

\usepackage{mathtools}
\usepackage{needspace}

\usepackage{titlesec}

\usepackage{latexsym}

\titleformat*{\section}{\large\bfseries}

\makeatletter
\renewcommand{\subsection}{\@startsection{subsection}{2}{0em}%
{\baselineskip}{-0em}{\bfseries\normalsize}}
\makeatother

\makeatletter
\newcommand{\listofdefinitions}{\@starttoc{def}}

\newcommand{\l@definition}[2]{\par\noindent#1 {\itshape #2}}
\makeatother

\newcounter{num}
\newenvironment{pflist}{\begin{list}{(\arabic{num})}{\usecounter{num} \leftmargin0cm \itemindent5pt}}{\end{list}}

\newcounter{num1}
\newenvironment{proplist}{\begin{list}{\arabic{num1}.}{\usecounter{num1} \listparindent0pt \topsep-2pt  \leftmargin25pt \itemindent0pt \itemsep5pt \partopsep0pt \labelwidth10pt \labelsep0.5em }}{\end{list}}

\newenvironment{parlist}{\begin{list}{(\arabic{num})}{\usecounter{num} \leftmargin0pt \itemindent5pt \topsep0pt \labelwidth0pt}}{\end{list}}

\theoremstyle{plain}
\theorembodyfont{\rmfamily}

\newtheorem{thm}{Theorem.}

\newtheorem{prop}[thm]{Proposition.}

\newtheorem{corollary}[thm]{Corollary.}

\newtheorem{pf}{Proof.}

\theoremstyle{break}
\theorembodyfont{\rmfamily}

\newcommand{\qed}{$\Box$}

\newcommand{\Aut}{\operatorname{Aut}}

\newcommand{\Span}{\operatorname{Span}}

\newcommand{\id}{\operatorname{id}}

\newcommand{\co}{\scriptstyle \operatorname{co}}

\newcommand{\deq}{\vcentcolon =}
\newcommand\coln{\colon}

\newcommand\1{{(1)}}
\newcommand\2{{(2)}}

\newcommand{\da}{\Delta_{A}}
\newcommand{\db}{\Delta_{B}}
\newcommand{\dhh}{\Delta_{H}}

\newcommand{\K}{1}
\newcommand{\ot}{\mathbin{\otimes}}
\newcommand{\ea}{{\varepsilon_{A}}}
\newcommand{\eb}{{\varepsilon_{B}}}
\newcommand{\eh}{{\varepsilon_{H}}}

\newcommand{\A}{1_A}
\newcommand{\B}{1_B}
\newcommand{\HH}{1_H}
\newcommand{\KK}{1_K}

\newcommand{\Z}{\mathbb Z}

\begin{document}

\maketitle

\begin{abstract}
\noindent
For two families of Yetter-Drinfel'd Hopf algebras considered earlier by the authors, we determine which of them are isomorphic. We also determine which of their Radford biproducts are isomorphic.
\end{abstract}

\section*{Introduction} \label{Sec:Introd}
\addcontentsline{toc}{section}{Introduction}
In a previous article~\cite{KaSo2}, the authors introduced two families of Yetter-Drinfel'd Hopf algebras over the Klein four-group, both of which depended on a not necessarily primitive fourth root of unity. Because the algebras in the first family were commutative, while the algebras in the second family were not commutative, no algebra in the first family can be isomorphic to an algebra in the second family. However, what was not addressed in the article was the question whether the Yetter-Drinfel'd Hopf algebras within one of these families are isomorphic, i.e., whether the choice of a different fourth root of unity actually leads to a different Yetter-Drinfel'd Hopf algebra. In this article, we show that this is actually the case; in other words, Yetter-Drinfel'd Hopf algebras arising from different fourth roots of unity are not isomorphic.

Every Yetter-Drinfel'd Hopf algebra gives rise to an ordinary Hopf algebra via the Radford biproduct construction. The arising biproducts were already described in our article~\cite{KaSo3}, where, in particular, we gave a presentation in terms of generators and relations and showed that the biproducts arising from the first family are isomorphic to the biproducts arising from the second family if in both cases the same root of unity is chosen. It is not easy to see whether or not biproducts corresponding to different roots of unity are isomorphic. In view of the result just mentioned, it suffices to discuss this question for the first family of Yetter-Drinfel'd Hopf algebras. As we show in Section~\ref{Sec:IsomRadf}, biproducts corresponding to two primitive fourth roots of unity are isomorphic, biproducts corresponding to two non-primitive fourth roots of unity are isomorphic, but two biproducts are not isomorphic if one of them corresponds to a primitive fourth root of unity, while the other corresponds to a non-primitive fourth root of unity. This is the main result of this article, which is stated in Theorem~\ref{IsomBiprodZetaOne}.

In the sequel, we assume that the reader is familiar with our previous work on the topic contained in the two references~\cite{KaSo2} and~\cite{KaSo3} already mentioned. The assumptions and the notation follow these references. We assume that our base field~$K$ is algebraically closed of characteristic zero, so that it contains in particular a primitive fourth root of unity, which is denoted by~$\iota$. All Yetter-Drinfel'd modules, and therefore all Yetter-Drinfel'd Hopf algebras, will be left-left Yetter-Drinfel'd modules; i.e., they are left modules and left comodules. The module action is denoted by a dot. The set of group-like elements of a Hopf algebra~$H$ is denoted by~$G(H)$. Its augmentation ideal, i.e., the kernel of the counit, is denoted by~$H^+$. The transpose of a linear map~$f$ is denoted by~$f^*$. Residue classes in a quotient space are denoted by a bar. The symbol~$\equiv$ is used to express that the residue classes of two vectors in a quotient space are equal. With respect to enumeration, we use the convention that propositions, definitions, and similar items are referenced by the paragraph in which they occur; i.e., a reference to Proposition~1.1 refers to the unique proposition in Paragraph~1.1.

Concerning the organization of the article, Section~\ref{Sec:IsomYet} covers the isomorphism problem for Yetter-Drinfel'd Hopf algebras. In Paragraph~\ref{YetDrinfFirst}, we briefly review how the first family of Yetter-Drinfel'd Hopf algebras, which is commutative, is defined. As stated, this definition depends on a fourth root of unity. Paragraph~\ref{IsomFirstZeta} compares two such algebras that correspond to a pair of fourth roots of unity that differ by a sign, leading up to the result in Paragraph~\ref{IsomFirstZetaOne} that Yetter-Drinfel'd Hopf algebras in the first family cannot be isomorphic if they correspond to two different roots of unity. In the case where the fourth root of unity under consideration is not primitive, we also derive there a certain equation of degree~$6$ for one of the generators. Paragraph~\ref{YetDrinfSec} briefly reviews the definition of the second family of Yetter-Drinfel'd Hopf algebras, which is noncommutative. As for the first family, two Yetter-Drinfel'd Hopf algebras in the second family cannot be isomorphic if they correspond to two different roots of unity, as shown in Paragraph~\ref{IsomSecZeta}. We also derive there an analogous equation of degree~$6$ for one of the generators, but now in the primitive case.

\enlargethispage{0.5mm}

Section~\ref{Sec:IsomRadf} treats the isomorphism problem for the corresponding Radford biproducts. Paragraph~\ref{RadfBiprod} describes Radford biproducts in general and lists some specific properties that are required for the proof of the main result. Paragraph~\ref{RadfBiprodCase} reviews the particular biproducts that arise from the Yetter-Drinfel'd Hopf algebras under consideration, which were first considered in~\cite{KaSo3}. Paragraph~\ref{PropRadfBiprod} and Paragraph~\ref{GroupLikeQuot} treat properties of these particular biproducts that will be needed in the sequel. Isomorphisms between different biproducts appear for the first time in Paragraph~\ref{IsomBiprodZeta}, where we construct an explicit isomorphism between the biproducts that arise from Yetter-Drinfel'd Hopf algebras corresponding to two fourth roots of unity that differ by a sign. Therefore, two Radford biproducts are isomorphic if the corresponding roots of unity are both primitive or both not primitive. The main result, namely that they are not isomorphic if one fourth root of unity is primitive and the other is not, is stated in Theorem~\ref{IsomBiprodZetaOne}. The rest of the article is concerned with the proof of that result. The fundamental technique of the proof is to analyze a certain group of group-like elements in the dual that has order~$4$. Under the transpose of a hypothetical isomorphism, it is mapped to another group of order~$4$, for which there are seven possibilities. These are treated on a case-by-case basis.

Both authors have presented results covered in this article, together with other results about the Yetter-Drinfel'd Hopf algebras considered here, at various conferences, either separately or jointly. The first author presented partial results at the AMS Spring Eastern Sectional Meeting in March~2021 and complete results at the AMS Spring Central Sectional Meeting in March~2022 as well as at the fourth `International Workshop on Groups, Rings, Lie and Hopf Algebras' in St.~John's in June~2022. The second author presented partial results at the Mathematical Congress of the Americas in July~2021 and complete results at the conference `Hopf Algebras, Monoidal Categories and Related Topics' in Bucharest in July~2022. Both authors presented together at the Joint Mathematics Meeting in April~2022 in two consecutive talks.

The travel of the second author to these conferences as well as his work on this article were supported by NSERC grant RGPIN-2017-06543. The work of the first author on this article was supported by a Faculty Summer Research Grant from the College of Science and Health at DePaul University.

\section{Isomorphisms of Yetter-Drinfel'd Hopf algebras} \label{Sec:IsomYet}
\subsection[The Yetter-Drinfel'd Hopf algebras in the first family]{} \label{YetDrinfFirst}
To begin, we recall the definition of the first family of examples of Yetter-Drinfel'd Hopf algebras that was introduced in~\cite{KaSo2}. If we use the notation \mbox{$\Z_2 = \{0,1\}$} for the group with two elements, the group
$G \deq \Z_2 \times \Z_2$ contains the four elements
\[g_1 \deq (0,0) \qquad g_2 \deq (1,0) \qquad g_3 \deq (0,1) \qquad  g_4 \deq (1,1)\]
with $g_1 = 1_G$ as its unit element if we write~$G$ multiplicatively, as we will in the sequel.
In \cite{KaSo2}, Sec.~2, the authors defined a commutative Yetter-Drinfel'd Hopf algebra~$A$ over the group algebra~$H \deq K[G]$ as the algebra generated by two commuting variables~$x$ and~$y$ subject to the defining relations
\[x^4 = 1 \qquad \qquad y^2 = \frac{1}{2}(1 + \zeta x + x^2 - \zeta x^3)\]
where~$\zeta$ is a not necessarily primitive fourth root of unity. We denote~$x$ by~$x_\zeta$ and~$y$ by~$y_\zeta$ if it is necessary to specify the root of unity under consideration, which happens in particular when two such algebras corresponding to different roots of unity need to be compared. In this situation, we will also denote~$A$ by~$A_\zeta$.

With the help of the primitive fourth root of unity~$\iota$ mentioned above, it is possible to introduce the elements $\omega_1 \deq 1$,
\[\omega_2 \deq \frac{\K}{2}(1 + \iota \zeta^2) x + \frac{\K}{2}(1 - \iota \zeta^2) x^3 \qquad
\omega_3 \deq \frac{\K}{2}(1 - \iota \zeta^2) x + \frac{\K}{2}(1 + \iota \zeta^2) x^3\]
$\omega_4 \deq x^2$, and
\[\eta_1 \deq y \qquad \eta_2 \deq x^3 y \qquad \eta_3 \deq x^2 y \qquad \eta_4 \deq x y\]
According to \cite{KaSo2}, Prop.~2.4, p.~105, these eight elements form a basis of~$A$, and the coproduct~$\da$ as well as the counit~$\ea$ are determined by the fact that these elements are group-like. As also shown there, the set $\{\omega_1 , \omega_2, \omega_3, \omega_4\}$ is closed under multiplication and in fact forms a group isomorphic
to~$\Z_2 \times \Z_2$. In~\cite{KaSo2}, Par.~2.2, p.~99, the 
$H$-action on~$A$ is introduced on the generators as
\[g_2.x = x^3 \qquad \qquad g_2.y = x^3 y \qquad \qquad
g_3.x = x \qquad \qquad g_3.y = x^2 y\]
which, according to~\cite{KaSo2}, Lem.~2.4, p.~107, means for the basis elements that
\[g_3.\omega_i = \omega_i \qquad \qquad  g_2.\omega_4 = \omega_4 \qquad \qquad g_2.\omega_2 = \omega_3 \qquad \qquad  g_i.\eta_{g_j} = \eta_{g_i g_j}\]
where we have used the notation~$\eta_{g_i} \deq \eta_i$. The coaction of~$H$ on~$A$ is derived from the action by defining
\begin{align*}
\delta(a) \deq \frac{1}{4} \sum_{g,g' \in G} \theta(g,g') \, g \ot g'.a 
\end{align*}
where~$\theta$ is the nondegenerate symmetric bicharacter described in~\cite{KaSo2}, Par.~2.2, p.~99 and determined by the condition
\[\begin{pmatrix}
\theta(g_2,g_2) & \theta(g_2,g_3) \\
\theta(g_3,g_2) & \theta(g_3,g_3)
\end{pmatrix} =
\begin{pmatrix}
\zeta^2 & -1 \\
-1 & 1
\end{pmatrix} \]
on its fundamental matrix.

\subsection[The first example for pairs of fourth roots of unity]{}
\label{IsomFirstZeta}
We begin by comparing the Yetter-Drinfel'd Hopf algebras~$A \deq A_\zeta$ and~$A' \deq A_{-\zeta}$ that arise from a fourth root of unity~$\zeta$ and its negative, regardless of whether~$\zeta$ is primitive or not. We will use the notation $x \deq x_\zeta$ and $y \deq y_\zeta$ for the generators of~$A$ and
$x' \deq x_{-\zeta}$ as well as $y' \deq y_{-\zeta}$ for the generators of~$A'$. These two algebras are not isomorphic:
\begin{prop} 
There is no isomorphism $f \coln A \to A'$ of Yetter-Drinfel'd Hopf algebras.
\end{prop}
\begin{pf} 
\begin{pflist} 
\item
Suppose that $f \coln A \to A'$ is such an isomorphism, which is \mbox{$H$-linear} and colinear by definition. We denote the corresponding basis consisting of group-like elements of~$A'$ by $\omega'_1$, $\omega'_2$, $\omega'_3$, $\omega'_4$, $\eta'_1$, $\eta'_2$, $\eta'_3$, and $\eta'_4$. Because $f$ preserves the unit element, we must have \mbox{$f(\omega_1) = \omega'_1$}, and since $\omega_4$ and $\omega'_4$ are, respectively, the only other group-like elements that are fixed points of the action, we must have $f(\omega_4) = \omega'_4$. Since~$f$ must map an orbit of length~$2$ among the group-like elements to an orbit of length~$2$, we must have $f(\omega_2) = \omega'_2$ or $f(\omega_2) = \omega'_3$. Because the action with~$g_2$ is an automorphism of~$A'$ that interchanges $\omega'_2$ and~$\omega'_3$, we can assume without loss of generality that $f(\omega_2) = \omega'_2$ and then also $f(\omega_3) = \omega'_3$. Because $\omega_1$, $\omega_2$, $\omega_3$, and~$\omega_4$ span the same space as $1$, $x$, $x^2$, and~$x^3$ (cf.~\cite{KaSo3}, Par.~2.1, p.~203), this implies that $f(x^k) = x'^k$ and in particular $f(x) = x'$.

\item
Because~$f$ maps group-like elements to group-like elements, we must have $f(\eta_1) = \eta'_i$ for some index~$i=0,1,2,3$, or equivalently $f(y) = x'^k y'$ for some number $k=0,1,2,3$. Since 
\begin{align*}
\frac{1}{2}(1_{A'} + \zeta x' + x'^2 - \zeta x'^3) = f(y^2) = x'^{2k} y'^2
= \frac{1}{2} x'^{2k} (1_{A'} - \zeta x' + x'^2 + \zeta x'^3)
\end{align*}
we see that we cannot have $k=0$ or~$k=2$, so $k=1$ or~$k=3$.

\item
Since we have
\begin{align*}
g_2.(x'y') = (g_2.x') (g_2.y') = x'^3 x'^3 y' = x'^2 y'
\end{align*}
we get in the case $k=1$, where $f(y) = x' y'$, that $g_2.f(y) = x'^2 y'$, but
\[f(g_2.y) = f(x^3 y) = f(x^3) f(y) = x'^3 x'y'= y'\]
so $k=1$ is not possible.

\item
Since $g_2.x'^2 = x'^2$, we get from the preceding computation that
\begin{align*}
g_2.(x'^3 y') = (g_2.x'^2) (g_2.(x'y')) = x'^2 x'^2 y' = y'
\end{align*}
In the case $k=3$, where $f(y) = x'^3 y'$, we therefore have $g_2.f(y) = y'$, but
\[f(g_2.y) = f(x^3 y) = f(x^3) f(y) =x'^3 x'^3 y'= x'^2 y'\] 
so $k=3$ is also not possible. Therefore, such an isomorphism~$f$ does not exist.~\qed
\end{pflist}
\end{pf}

\subsection[Primitive and non-primitive fourth roots of unity in the first family]{}
\label{IsomFirstZetaOne}
Suppose now that~$\zeta$ is a primitive fourth root of unity, and denote by \mbox{$A'' \deq A_\xi$} the analogous algebra, defined with a non-primitive fourth root of unity~$\xi = \pm 1$. Then these two algebras are not isomorphic as Yetter-Drinfel'd Hopf algebras: Under the action of the group, both contain precisely one orbit of length~$4$ consisting of group-like elements. But it follows from~\cite{KaSo2}, Lem.~2.1, p.~96 that the group-like elements in this orbit have order exactly~$4$ in the case of~$A$ and order exactly~$8$ in the case of~$A''$. This completes the relatively simple~proof of our first main result about the first family of examples of Yetter-Drinfel'd Hopf algebras:
\begin{thm}
If~$\zeta$ and~$\xi$ are two distinct fourth roots of unity, the Yetter-Drinfel'd Hopf algebras~$A_\zeta$ and~$A_\xi$ are not isomorphic.
\end{thm}

The preceding argument relied on~\cite{KaSo2}, Lem.~2.1, which states in particular that~$y^4 = x^2$ if~$\zeta^2 = 1$, so that~$y^8 = 1$. 
It is worth noting that we have in fact the following stronger relation:
\begin{prop}
If~$\zeta = \pm 1$, we have
\[y^6 - y^4 + y^2 - 1 = 0\]
\end{prop}
\begin{pf}
As just mentioned, we have $y^4 = x^2$. As a consequence of the defining relation
$y^2 = \frac{\K}{2}(1 + \zeta x + x^2 - \zeta x^3)$, we therefore have
\begin{align*}
y^6 = \frac{\K}{2}(1 - \zeta x + x^2 + \zeta x^3)
\end{align*}
which implies that
$y^2 + y^6 = 1 + x^2 = 1 + y^4$.
\qed
\end{pf}

This relation is stronger in the following sense: The equation~$y^8=1$ means that~$y^2$ satisfies the polynomial
\[X^4-1 = (X^2+1) (X^2-1) = (X^2+1) (X+1) (X-1)\]
in the indeterminate~$X$. The proposition states that~$y^2$ already satisfies its proper divisor
\[X^3 - X^2 + X - 1 = (X^2+1) (X-1)\]
which has one cyclotomic polynomial less in its factorization (cf.~\cite{Jac3}, Sec.~III.1, p.~111).

\subsection[The Yetter-Drinfel'd Hopf algebras in the second family]{} \label{YetDrinfSec}
We now turn to a second, similar family of Yetter-Drinfel'd Hopf algebras, which was introduced in~\cite{KaSo2}, Sec.~3. These algebras are also defined over the group algebra~$H \deq K[G]$ of the Klein four-group, but are no longer commutative. The~algebras in this second family are also distinguished by their dependence on a not necessarily primitive fourth root of unity~$\zeta$, and will also be denoted by~$A_\zeta$, or by~$A$ if a fixed fourth root of unity is considered. By definition,~$A$~is generated by two noncommuting variables~$x$ and~$y$ subject to the defining relations
\[x^4 = 1 \qquad \qquad xy = yx^3 \qquad \qquad 
y^2 = \frac{\K}{2}(\zeta 1 + x - \zeta x^2 + x^3)\]
As for the first family, we denote~$x$ by~$x_\zeta$ and~$y$ by~$y_\zeta$ if it is necessary to specify the root of unity under consideration. The elements $\omega_1$, $\omega_2$, $\omega_3$, $\omega_4$, $\eta_1$, $\eta_2$, $\eta_3$, and~$\eta_4$ are defined by the same formulas as in Paragraph~\ref{YetDrinfFirst} and form again a basis of~$A$, as shown in \cite{KaSo2}, Prop.~3.4, p.~115. Again, the coalgebra structure is determined by the fact that these elements are group-like.

The $H$-action on~$A$ is constructed in a similar way as the action for the first family. It is given on the generators by the formulas
\[g_2.x = x^3 \qquad \qquad g_2.y = x^3 y \qquad \qquad
g_3.x = x \qquad \qquad g_3.y = x^2 y\]
(cf.~\cite{KaSo2}, Par.~3.2, p.~111). As also explained there, the $H$-coaction is constructed from the $H$-action in exactly the same way as in the case of the first family reviewed in Paragraph~\ref{YetDrinfFirst}, and the $H$-action permutes the group-like elements also as in the case of the first family, as recorded in~\cite{KaSo2}, Lem.~3.2, p.~117.

\subsection[The second example for pairs of fourth roots of unity]{}
\label{IsomSecZeta}
The proof of the analogue of Theorem~\ref{IsomFirstZetaOne} for the second family of Yetter-Drinfel'd Hopf algebras is easier:
\begin{thm}
If~$\zeta$ and~$\xi$ are two distinct fourth roots of unity, the Yetter-Drinfel'd Hopf algebras~$A_\zeta$ and~$A_\xi$ are not isomorphic.
\end{thm}
\begin{pf} 
We use the notations introduced in Paragraph~\ref{YetDrinfSec} for~$A \deq A_\zeta$ and use primes for the corresponding elements in~$A' \deq A_\xi$, so that in particular~$x' \deq x_\xi$ and~$y' \deq y_\xi$. Suppose that $f \coln A \to A'$ is an isomorphism. The same argument used in the proof of Proposition~\ref{IsomFirstZeta} yields that we can assume that~$f(x) = x'$. Because~$f$ maps group-like elements to group-like elements, we have as in the previous case that $f(\eta_1) = \eta'_i$ for some index~$i=0,1,2,3$, or equivalently $f(y) = x'^k y'$ for some number $k=0,1,2,3$. But since 
\begin{align*}
\frac{1}{2}(\zeta 1_{A'} +  x' - \zeta x'^2 + x'^3) &= f(y^2) = f(y)^2 = x'^{k} y' x'^k y' 
= y' x'^{3k} x'^k y' = y'^2\\
&= \frac{1}{2} (\xi 1_{A'} + x' - \xi x'^2 + x'^3)
\end{align*}
this implies that~$\zeta = \xi$, contrary to our assumption.~\qed
\end{pf}

If~$\zeta$ is a primitive fourth root of unity and~$\xi = \pm 1$ is not, one can alternatively argue as in Paragraph~\ref{IsomFirstZetaOne}: Under the action of the group, both~$A$ and~$A'$ contain precisely one orbit of length~$4$ consisting of group-like elements. We have just seen in the preceding proof that these group-like elements have all the same square, so it follows from~\cite{KaSo2}, Lem.~3.1, p.~109 that the group-like elements in this orbit have order exactly~$8$ in the case of~$A$ and order exactly~$4$ in the case of~$A'$. Note that in this lemma, the roles of the primitive and the non-primitive fourth roots of unity are reversed in comparison to~\cite{KaSo2}, Lem.~2.1: For the second family, we have~$y^4 = x^2$ if~$\zeta$ is primitive, so that~$y^8 = 1$. Again, there is a stronger relation:
\begin{prop}
If~$\zeta$ is primitive, we have
\[y^6 - \zeta y^4 - y^2 + \zeta 1 = 0\]
\end{prop}
\begin{pf}
As just mentioned, we have $y^4 = x^2$. As a consequence of the defining relation
$y^2 = \frac{\K}{2}(\zeta 1 + x - \zeta x^2 + x^3)$, we therefore have
\begin{align*}
y^6 = \frac{\K}{2}(- \zeta 1 + x + \zeta x^2 + x^3)
\end{align*}
which implies that
$y^6 - y^2 = - \zeta 1 + \zeta x^2 = - \zeta 1 + \zeta y^4$.
\qed
\end{pf}

This relation is stronger in the following sense: The equation~$y^8=1$ means that~$y^2$ satisfies the polynomial
\[X^4-1 = (X - \zeta) (X - \zeta^2)  (X - \zeta^3) (X - \zeta^4) 
= (X - \zeta) (X + 1)  (X + \zeta) (X - 1)\]
in the indeterminate~$X$. The proposition states that~$y^2$ already satisfies its proper divisor
\[X^3 - \zeta X^2 - X + \zeta = (X - \zeta) (X^2 - 1) = (X - \zeta) (X + 1) (X - 1)\]
which has one factor less.

\section{Isomorphisms of Radford biproducts} \label{Sec:IsomRadf}
\subsection[Radford biproducts]{} \label{RadfBiprod}
Associated with any Yetter-Drinfel'd Hopf algebra~$A$ over a Hopf algebra~$H$ is an ordinary Hopf algebra~$B \deq A \star H$, which was introduced in~\cite{RadfProj} and is called the Radford biproduct of~$A$ and~$H$. It has the tensor product~$A \ot H$ as underlying vector space, but it
is endowed with the so-called smash product multiplication 
\begin{align*}
( a \star h) ( a' \star  h' ) \deq
a ( h_{\1}.a') \star h_{\2} h'
\end{align*}
which is in general different from the usual tensor product multiplication. The comultiplication is the so-called cosmash product comultiplication
\begin{align*} 
\db ( a \star h) \deq
(a_{\1} \star a_{\2}{}^\1 h_{\1} ) \ot (a_{\2}{}^\2 \star h_{\2} )
\end{align*}
Here we have used Heyneman-Sweedler sigma notation in the form 
\[\da(a) = a_{\1} \ot a_{\2} \in A \ot A \qquad \qquad
\delta(a) = a^{\1} \ot a^{\2} \in H \ot A\]
with lower indices for the coproduct~$\da$ of~$A$, and also for the coproduct~$\dhh$ of~$H$, and with upper indices for the coaction~$\delta$. Further details and references are given in~\cite{KaSo3}, Par.~1.2, p.~198, where also explicit formulas for unit, counit, and antipode are stated. Note that the map
$H \to B,~h \mapsto \A \ot h$
is an injective Hopf algebra homomorphism, and we will at times view~$H$ as a Hopf subalgebra of~$B$ via this map.

The original algebra~$A$ is in fact a left module over its Radford biproduct~$B$. This can be understood from the point of view of an induced module construction: We consider the base field~$K$ as a trivial $H$-module via the counit~$\eh$ of~$H$ and form the induced module~$B \ot_H K$. Since
\[B \ot_H K = (A \star H) \ot_H K \cong A\]
this induced module is isomorphic to~$A$, and we can pull back the $B$-module structure along this isomorphism. As we have
\begin{align*}
(a \star h&)(a' \star 1_H) \ot_H 1_K = (a(h_\1.a') \star h_\2) \ot_H 1_K \\
&= (a(h_\1.a') \star \HH) \ot_H \eh(h_\2) \KK 
= (a(h.a') \star \HH) \ot_H \KK
\end{align*}
the arising $B$-action on~$A$ is given by
\[(a \star h).a' = a(h.a')\]
In particular, the restriction of the $B$-action on~$A$ to~$H$ agrees with the original action in the sense that $(\A \star h).a = h.a$. With respect to this action, $A$ is a left $B$-module coalgebra, as it is easy to see that the counit~$\ea$ of~$A$ is $B$-linear, and the coproduct~$\da$ of~$A$ is $B$-linear because
\begin{align*}
\da((a \star h).a') &= \da(a(h.a')) = a_\1 ((a_\2{}^\1 h_\1).a'_\1) \ot (a_\2{}^\2 (h_\2.a'_\2)) \\
&= ((a_\1 \star a_\2{}^\1 h_\1).a'_\1) \ot ((a_\2{}^\2 \star h_\2).a'_\2) \\
&= ((a \star h)_\1.a'_\1) \ot ((a \star h)_\2.a'_\2)
\end{align*}

The spaces~$B \ot_H K$ and~$A$ can be compared to the space~$B/BH^+$, which is clearly a left $B$-module, and also a quotient coalgebra. It is easy to see that its coproduct is $B$-linear, so that $B/BH^+$ is also a left $B$-module coalgebra. 
\begin{prop}
The map 
\[A \to B/BH^+,~a \mapsto \overline{a \star \HH}\]
is an isomorphism of left $B$-module coalgebras.
\end{prop}
\begin{pf}
Note that $BH^+ = (A \star H)H^+ = A \star (HH^+) = A \star H^+$. This shows that the kernel of the map 
\[B\to A,~a \star h \mapsto \eh(h) a\]
is precisely $BH^+$, and the induced map from the quotient space~$B/BH^+$ to~$A$
is the inverse of the map under consideration,
as we have $\overline{a \star h} = \overline{a \star \eh(h) \HH}$ in the quotient $B/BH^+$. Now we get
\begin{align*}
\Delta_{B/BH^+}\left(\overline{a \star \HH}\right) 
&= \overline{(a \star \HH)_\1} \ot \overline{(a \star \HH)_\2}
= \overline{(a_\1 \star a_\2{}^\1)} \ot \overline{(a_\2{}^\2 \star \HH)} \\
&= \overline{(a_\1 \star \eh(a_\2{}^\1) \HH)} \ot \overline{(a_\2{}^\2 \star \HH)} \\
&= \overline{(a_\1 \star \HH)} \ot \overline{(a_\2 \star \HH)}
\end{align*}
which proves that our map is a coalgebra isomorphism, as it is clearly also compatible with the counits. Since 
\begin{align*}
\overline{((a \star h).a') \star \HH} &= \overline{a(h.a') \star \HH} 
= \overline{a(h_\1.a') \star h_\2} \\
&= \overline{(a \star h) (a' \star \HH)} = (a \star h) \overline{(a' \star \HH)}
\end{align*}
this isomorphism is also $B$-linear.~\qed
\end{pf}

With the help of the map
\[\pi_H \coln B \to H,~a \star h \mapsto \ea(a) h\]
we can define the space
$B^{\co H} \deq \{b \in B \mid (\id_B \ot \pi_H)(\db(b)) = b \ot \HH \}$
of coinvariants. The preceding proposition yields a bijection between the space of coinvariants and the above quotient space:
\begin{corollary}
The map 
\[B^{\co H} \to B/BH^+,~b \mapsto \overline{b}\]
is bijective.
\end{corollary}
\begin{pf}
Note that
\begin{align*}
(\id_B \ot \pi_H)(\db (a \star h)) &=
(a_{\1} \star a_{\2}{}^\1 h_{\1} ) \ot \ea(a_{\2}{}^\2) h_{\2} = (a \star h_{\1} ) \ot  h_{\2} 
\end{align*}
From this fact, it follows easily that the map
\[A \to B^{\co H},~a \mapsto a \star \HH\]
is bijective (cf.~\cite{BM}, Lem.~1.5, p.~39; \cite{HS}, Lem.~2.5.2, p.~89), and so the preceding proposition implies the assertion.
\qed
\end{pf}

\subsection[The Radford biproduct in our case]{} \label{RadfBiprodCase}
The Radford biproduct construction can be applied in particular to the two families of Yetter-Drinfel'd Hopf algebras considered in Section~\ref{Sec:IsomYet}, where \mbox{$H=K[G]$} was the group algebra of the Klein four-group~$G$, and~$A$ depended on a fourth root of unity~$\zeta$ and was therefore denoted in both cases by~$A_\zeta$. Accordingly, we will denote the corresponding Radford biproduct~$B$ also by~$B_\zeta$ if we want to emphasize the dependence on~$\zeta$.

It was shown in~\cite{KaSo3}, Cor.~7.2, p.~222 that the Radford biproduct~$B_\zeta$ arising from the first family is isomorphic to the Radford biproduct arising from the second family that corresponds to the same fourth root of unity~$\zeta$. To determine which of these Radford biproducts are isomorphic, it is therefore sufficient to focus on the first family. The basic properties of such a Radford biproduct were described in~\cite{KaSo3}, Sec.~2, p.~202ff: It is generated by the four elements
\[u \deq x \star \HH \qquad \quad v \deq y \star \HH \qquad \quad
r \deq \A \star g_2 \qquad \quad s \deq \A \star g_3\]
\pagebreak

which satisfy the relations
\begin{proplist}
\item
$\displaystyle
u^4 = 1, \qquad uv = vu, \qquad v^2 = \frac{1}{2}(1 + \zeta u + u^2 - \zeta u^3)$

\item
$\displaystyle
r^2 = 1, \qquad rs = sr, \qquad s^2 = 1$

\item
$\displaystyle
r u = u^3 r, \qquad r v = u^3 v r, \qquad s u = u s, \qquad s v = u^2 v s$
\end{proplist}
and these relations are defining. The elements $u^i v^j r^k s^l$, where the index~$i$ takes the values~$0$, $1$, $2$, and~$3$ and the other three indices $j$, $k$, and~$l$ take the values~$0$ and~$1$, form a basis of~$B$. We will denote~$u$, $v$, $r$, and~$s$ by~$u_\zeta$, $v_\zeta$, $r_\zeta$, and~$s_\zeta$ if we want to emphasize that they are the generators of~$B_\zeta$. 

The formula for the cosmash product comultiplication directly yields the values of the coproduct on the generators, which are stated in~\cite{KaSo3}, Prop.~2.3, p.~204, together with the values of the counit on the generators. The generators~$r$ and~$s$ are group-like, but the generators~$u$ and~$v$ are not. Rather, we have
\begin{align*}
\db(u) &= \frac{\K}{2} (u \ot u + u \ot u^3 +  u^3 s \ot u - u^3 s \ot u^3)
\end{align*}
and
\begin{align*}
\db(v) &= \frac{1}{4} v (1 + r + s + rs) \ot v
+ \frac{1}{4} v (1 - \zeta^2 r - s + \zeta^2 rs) \ot u v \\
&\quad + \frac{1}{4} v (1 - r + s - rs) \ot u^2 v
+ \frac{1}{4} v (1 + \zeta^2 r - s - \zeta^2 rs) \ot u^3 v
\end{align*}
The counit is given on generators by~$\eb(u) = \eb(v) = \eb(r) = \eb(s) = 1$. Formulas for the values of the antipode on the generators are given in~\cite{KaSo3}, Prop.~2.4, p.~204.

The group-like elements in~$A$ and~$H$ lead to corresponding elements of~$B$. As in~\cite{KaSo3}, Par.~2.2, p.~203, we denote them by
\[c_i \deq \omega_i \star \HH \qquad  \qquad  d_i \deq \eta _i \star \HH \qquad \qquad
h_i \deq 1_A \star g_i\]
for $i=1,2,3,4$.

\subsection[Properties of the Radford biproduct]{} \label{PropRadfBiprod}
For any finite-dimensional Hopf algebra~$B$, an algebra homomorphism $\chi \in G(B^*)$ to the base field~$K$ yields an algebra automorphism 
\[\psi_\chi \coln B \to B,~b \mapsto \psi_\chi(b) \deq (\id_B \ot \chi)(\db(b)) 
= b_\1 \chi(b_\2)\]
(cf.~\cite{SoRib}, Prop.~2.6, p.~427), and the assignment
\[G(B^*) \to \Aut(B),~\chi \mapsto \psi_\chi\]
is a group homomorphism. If~$B$ is the specific Radford biproduct just described in Paragraph~\ref{RadfBiprodCase}, the group~$G(B^*)$ has been determined in~\cite{KaSo3}, Prop.~3.1, p.~207. It is elementary abelian of order~$8$ and generated by three specific elements~$\chi_1$,~$\chi_2$, and~$\chi_3$ given there. We use the abbreviated notation
$\psi_i \deq \psi_{\chi_i}$ for~$i=1,2,3$. The basis elements introduced in Paragraph~\ref{RadfBiprodCase} are eigenvectors of these automorphisms:
\begin{prop} 
Suppose that $b = u^i v^j r^k s^l$, where the index~$i$ is~$0$, $1$, $2$, or~$3$ and the other three indices $j$, $k$, and~$l$ are~$0$ or~$1$. Then we have
\begin{enumerate}
\item
$\psi_1(b) = (-1)^j b$

\item 
$\psi_2(b) = (-1)^k b$

\item 
$\psi_3(b) = (-1)^l b$
\end{enumerate}
\end{prop}
\begin{pf} 
\begin{pflist} 
\item
For the first statement, the formula for the coproduct of~$u$ recalled in Paragraph~\ref{RadfBiprodCase} yields
\begin{align*}
\psi_1(u) &= \frac{1}{2} (u \chi_1(u) + u \chi_1(u^3) +  u^3 s \chi_1(u) - u^3 s \chi_1(u^3)) 
= u
\end{align*}
as we have $\chi_1(u) = 1$. Similarly, the formula for the coproduct of~$v$ also recalled there yields
\begin{align*}
\psi_1(v) &= \frac{1}{4} v (1 + r + s + rs) \chi_1(v)
+ \frac{1}{4} v (1 - \zeta^2 r - s + \zeta^2 rs) \chi_1(u v) \\
&\quad + \frac{1}{4} v (1 - r + s - rs) \chi_1(u^2 v)
+ \frac{1}{4} v (1 + \zeta^2 r - s - \zeta^2 rs) \chi_1(u^3 v)
\end{align*}
Since $\chi_1(v) = \chi_1(u v) = \chi_1(u^2 v) = \chi_1(u^3 v) = -1$, this expression reduces to~$\psi_1(v) = -v$. As it is immediate that~$\psi_1(r) = r$ and~$\psi_1(s) = s$, the first assertion holds.

\item
From the construction of the characters~$\chi_2$ and~$\chi_3$ in~\cite{KaSo3}, Prop.~3.1, p.~207, it follows that they are pullbacks of characters in~$G(H^*)$ along the map~$\pi_H$ considered in Paragraph~\ref{RadfBiprod}. Since~$u$ and~$v$ are contained in $B^{\co H}$, this implies that $\psi_2(u) = \psi_3(u) = u$ and $\psi_2(v) = \psi_3(v) =v$. For the two group-like elements~$r$ and~$s$, we have $\psi_2(r) = -r$, $\psi_2(s) = s$, $\psi_3(r) = r$, and~$\psi_3(s) = -s$. This implies the second and the third assertion. 
\qed
\end{pflist} 
\end{pf} 

The preceding proposition implies that the span of the powers of~$u$ consists exactly of those elements that are invariant under the above action of the group-like elements of~$B^*$: 
\begin{corollary} 
\[\Span(1, u, u^2, u^3) = \{b \in B \mid \psi_\chi(b) = b \; \text{for all} \; \chi \in G(B^*)\}\]
\end{corollary}
\begin{pf} 
Because the assignment $\chi \mapsto \psi_\chi$ is a group homomorphism, an element~$b \in B$ is contained in the set on the right-hand side if and only if 
\[\psi_1(b) =  \psi_2(b) = \psi_3(b) = b\]
Clearly, it follows from the preceding proposition that the powers of~$u$ have this property. But by expanding a given element~$b \in B$ in terms of the basis elements~$u^i v^j r^k s^l$ considered in that proposition, it follows that~$b$ can only have this property if the coefficients of the basis elements with~$j \neq 0$, $k \neq 0$, or~$l \neq 0$ vanish. But this means that~$b$ is a linear combination of the elements~$u^i$.
\qed
\end{pf}

More generally, the above proposition implies that~$B$ is the direct sum of the simultaneous eigenspaces of all the mappings~$\psi_\chi$ for~$\chi \in G(B^*)$, and that these eigenspaces are spanned by those basis elements~$u^i v^j r^k s^l$ that are contained in that eigenspace.

\subsection[Group-like elements in quotients]{} \label{GroupLikeQuot}
According to~\cite{KaSo3}, Lem.~3.1, p.~207, the group of group-like elements~$G(B)$ of~$B$ is generated by~$u^2$, $r$, and~$s$. The subgroup generated by~$r$ and~$s$ spans the subalgebra of~$B$ that corresponds to~$H$. As we have discussed in Paragraph~\ref{RadfBiprod}, the quotient~$B/BH^+$ is isomorphic to~$A$ as a left $B$-module coalgebra, and under this isomorphism, $y \in A$ corresponds to
$\overline{v} \in B/BH^+$, so $\overline{v}$ is group-like, and clearly $\overline{vs} = \overline{v}$, so the residue class of~$vs$ is group-like, too.

The subgroup $\langle u^2 r, s \rangle$ of~$G(B)$ spans a Hopf subalgebra~$\tilde{H} \deq \Span(\langle u^2 r, s \rangle)$ to which a similar statement applies: 
\begin{prop}
In the quotient~$B/B \tilde{H}^+$, the residue class 
\[\frac{1+\iota}{2} \overline{v} + \frac{1-\iota}{2} \overline{u^2 v} =
\frac{1+\iota}{2} \overline{v s} + \frac{1-\iota}{2}  \overline{u^2 v s}\]
is group-like.
\end{prop}
\begin{pf}
As recalled in Paragraph~\ref{RadfBiprodCase}, we have
\begin{align*}
\db(v) &= \frac{1}{4} v (1 + r + s + rs) \ot v
+ \frac{1}{4} v (1 - \zeta^2 r - s + \zeta^2 rs) \ot u v \\
&\quad + \frac{1}{4} v (1 - r + s - rs) \ot u^2 v
+ \frac{1}{4} v (1 + \zeta^2 r - s - \zeta^2 rs) \ot u^3 v
\end{align*}
Now we have $v (s-1) \in B\tilde{H}{}^+$ and therefore $v \equiv v s \pmod{B\tilde{H}{}^+}$. Similarly, $v u^2 (u^2 r-1) \in B\tilde{H}{}^+$ and 
$v u^2 (u^2 rs -1) \in B\tilde{H}{}^+$ and therefore 
\[v r \equiv v rs \equiv u^2 v \pmod{B\tilde{H}{}^+}\]
For the quotient coalgebra $B/B\tilde{H}{}^+$, this means that
\begin{align*}
\Delta(\overline{v}) &= \frac{1}{2} (\overline{v} + \overline{u^2 v}) \ot \overline{v}  + \frac{1}{2} (\overline{v} - \overline{u^2 v}) \ot \overline{u^2 v} \\
&= \frac{1}{2} (\overline{v}  \ot \overline{v}  
+ \overline{u^2 v} \ot \overline{v}  
+ \overline{v} \ot \overline{u^2 v}
- \overline{u^2 v} \ot \overline{u^2 v})
\end{align*}
Since $B/B\tilde{H}{}^+$ is a left $B$-module coalgebra, we have
\[\Delta(\overline{u^2 b}) = \Delta(u^2.\overline{b}) = \Delta(u^2).\Delta(\overline{b}) = 
(u^2 \ot u^2).\Delta(\overline{b})\] 
so the last formula also yields that 
\begin{align*}
\Delta(\overline{u^2 v}) 
&= \frac{1}{2} (-\overline{v}  \ot \overline{v}  
+ \overline{u^2 v} \ot \overline{v}  
+ \overline{v} \ot \overline{u^2 v}
+ \overline{u^2 v} \ot \overline{u^2 v})
\end{align*}

This implies 
\begin{align*}
\Delta \left(\frac{1+\iota}{2} \overline{v} 
+ \frac{1-\iota}{2}\overline{u^2 v} \right) 
&= \frac{1}{2} \left(\iota \overline{v}  \ot \overline{v}  
+ \overline{u^2 v} \ot \overline{v}  
+ \overline{v} \ot \overline{u^2 v}
- \iota \overline{u^2 v} \ot \overline{u^2 v}\right) \\
&= \left(\frac{1+\iota}{2} \overline{v} + \frac{1-\iota}{2}\overline{u^2 v}\right)  \ot
\left(\frac{1+\iota}{2} \overline{v} + \frac{1-\iota}{2}\overline{u^2 v}\right) 
\end{align*}
So the residue class of~$\frac{1+\iota}{2} v + \frac{1-\iota}{2} u^2 v$, which is equal to the residue class of~$\frac{1+\iota}{2} v s + \frac{1-\iota}{2} u^2 v s$, is group-like in $B/B\tilde{H}{}^+$.
\qed
\end{pf}

\subsection[The biproducts for pairs of fourth roots of unity]{}
\label{IsomBiprodZeta}
Compared to the Yetter-Drinfel'd Hopf algebras in Section~\ref{Sec:IsomYet}, it is substantially more complicated to determine which of their biproducts are isomorphic. However, it is still relatively easy to derive a variant of Proposition~\ref{IsomFirstZeta} for the biproducts: The biproducts~$B \deq B_\zeta$ and~$B' \deq B_{-\zeta}$ are isomorphic, regardless of whether~$\zeta$ is primitive or not. We will use the notation $u \deq u_\zeta$ and $u' \deq u_{-\zeta}$, and analogous notations for the other generators~$v$, $r$, and~$s$.
\begin{prop} 
There is a Hopf algebra isomorphism $f \coln B \to B'$ that takes the values
\begin{align*}
f(u) = u'^3 \qquad f(v) = v' \qquad f(r) = r's' \qquad f(s) = s' 
\end{align*}
on the generators.
\end{prop}
\begin{pf} 
\begin{pflist} 
\item
To establish the existence of an algebra homomorphism~$f$ that takes the proposed values on the generators, we have to check that these proposed values satisfy the defining relations of~$B$ stated in Paragraph~\ref{RadfBiprodCase}. In the first set of relations, this is obvious for the first two, and for the third we have
\begin{align*}
f(v)^2 = v'^2 = \frac{1}{2}(1 - \zeta u' + u'^2 + \zeta u'^3)
= \frac{1}{2}(1 + \zeta f(u) + f(u)^2 - \zeta f(u)^3)
\end{align*}
It is obvious that the second set of relations is satisfied, and for the third set, we have
\begin{align*}
f(r) f(u) = r' s' u'^3 = u' r' s'  = f(u)^3 f(r)
\end{align*}
for the first relation and
\begin{align*}
f(r) f(v) &= r' s' v' = r' u'^2 v' s' =  u'^2 r' v' s' =  u'^2 u'^3 v' r' s' \\
&= u' v' r' s' = (u'^3)^3 v' r' s' = f(u)^3 f(v) f(r)
\end{align*}
for the second relation. The third relation is obvious, and for the fourth relation, we have
$f(s) f(v) = s' v' = u'^2 v' s' = f(u)^2 f(v) f(s)$.

\item
To see that~$f$ is also a coalgebra homomorphism, we need to establish the 
equations~$\Delta_{B'}(f(b)) = (f \ot f)(\db(b))$ and~$\varepsilon_{B'}(f(b)) = \eb(b)$. Because in both equations both sides depend multiplicatively on~$b$, it suffices to verify them when~$b$ is one of the generators, using the formulas recalled in Paragraph~\ref{RadfBiprodCase}. In the case of the counit, this verification is immediate. In the case of the coproduct, this verification is also immediate for the two generators~$r$ and~$s$, which are group-like. Using that~$u'^2$ is also a group-like element, as already mentioned in Paragraph~\ref{GroupLikeQuot}, we have
\begin{align*}
\Delta_{B'}(f(u)) &= \Delta_{B'}(u'^3) = (u'^2 \ot u'^2) \Delta_{B'}(u') \\
&= \frac{\K}{2} (u'^2 \ot u'^2) (u' \ot u' + u' \ot u'^3 +  u'^3 s' \ot u' - u'^3 s' \ot u'^3)\\
&= \frac{\K}{2} (u'^3 \ot u'^3 + u'^3 \ot u' +  u' s' \ot u'^3 - u' s' \ot u') 
\end{align*}
Since
\begin{align*}
(f \ot f) (\db(u)) 
&= \frac{\K}{2} (f \ot f) (u \ot u + u \ot u^3 +  u^3 s \ot u - u^3 s \ot u^3) 
\end{align*}
this shows that $\Delta_{B'}(f(u)) = (f \ot f) (\db(u))$. 

\item
On the right-hand side of the equation
\begin{align*}
\db(v) &= \frac{1}{4} (v \ot v) \bigl((1 + r + s + rs) \ot 1
+  (1 - \zeta^2 r - s + \zeta^2 rs) \ot u  \\
&\quad + (1 - r + s - rs) \ot u^2 
+ (1 + \zeta^2 r - s - \zeta^2 rs) \ot u^3\bigr)
\end{align*}
recalled in Paragraph~\ref{RadfBiprodCase}, the last factor is a sum of four terms. A similar equation holds for the coproduct of~$v'$. In a sense, $f \ot f$ preserves the first and the third term, as we have
\[f(1 + r + s + rs) = 1 + r' + s' + r's' \quad  \quad 
f(1 - r + s - rs) = 1 - r' + s' - r's'\]
and interchanges the second and the fourth term, as we have
\[f(1 - \zeta^2 r - s + \zeta^2 rs) = 1 + \zeta^2 r' - s' - \zeta^2 r's'\]
and~$f(1 + \zeta^2 r - s - \zeta^2 rs) = 1 - \zeta^2 r' - s' + \zeta^2 r's'$. This implies that
\[(f \ot f) (\db(v)) = \Delta_{B'}(v') = \Delta_{B'}(f(v)) \]
completing the proof that~$f$ is also a coalgebra homomorphism. Since~$f$ then automatically commutes with the antipode (cf.~\cite{Sweed}, Lem.~4.0.4, p.~81), it is also a Hopf algebra homomorphism. By construction, the image of~$f$ contains a set of generators of~$B'$, so~$f$ is surjective, and thus bijective for dimension considerations.~\qed
\end{pflist}
\end{pf}

There is another, slightly more explicit way of showing that~$f$ is bijective: By interchanging the roles of~$\zeta$ and~$-\zeta$ in the preceding proposition, we obtain a Hopf algebra homomorphism
$f' \coln B' \to B$ that takes the values
\begin{align*}
f'(u') = u^3 \qquad f'(v') = v \qquad f'(r') = rs \qquad f'(s') = s 
\end{align*}
on the generators. Then~$f' \circ f$ preserves the generators~$u$, $v$, $r$, and~$s$ and is therefore the identity mapping. Similarly, $f \circ f'$ is the identity mapping, which implies that~$f$ is bijective with inverse~$f'$.

\subsection[Primitive and non-primitive fourth roots of unity for biproducts]{}
\label{IsomBiprodZetaOne}
If~$\zeta$ and~$\xi$ are two fourth roots of unity, it follows from Proposition~\ref{IsomBiprodZeta} that~$B_\zeta$ and~$B_\xi$ are isomorphic if both~$\zeta$ and~$\xi$ are primitive and that they are also isomorphic if both~$\zeta$ and~$\xi$ are not primitive. The converse of this result is also true:
\begin{thm}
If~$\zeta$ and~$\xi$ are two fourth roots of unity, the Radford biproducts~$B_\zeta$ and~$B_\xi$ are isomorphic if and only if both~$\zeta$ and~$\xi$ are primitive or both~$\zeta$ and~$\xi$ are not primitive.
\end{thm}

The rest of this article is concerned with the proof of this theorem. So, for two fourth roots of unity~$\zeta$ and~$\xi$, we suppose that $f \coln B \deq B_\zeta \to B' \deq B_\xi$ is an isomorphism of Hopf algebras. We will often use primes instead of the index~$\xi$ to refer to the quantities corresponding to $B_\xi$. First, we determine the possible images of~$u_\zeta$ and its powers: 
\begin{prop} 
We have $f(u_\zeta^2) = u_\xi^2$. Furthermore, we have 
$f(u_\zeta) = u_\xi$ or $f(u_\zeta) = u_\xi^3$.
\end{prop}
\begin{pf}
\begin{pflist}
\item
By \cite{KaSo3}, Lem.~3.1, p.~207, $u_\zeta^2$ is the only nontrivial central group-like element in~$B$, so it must be mapped to the only nontrivial central group-like element of~$B'$. This proves the first assertion.

\item
Since $f^*(G(B_\xi^*)) = G(B_\zeta^*)$, the description in Corollary~\ref{PropRadfBiprod} shows that $f(\Span(1, u_\zeta, u_\zeta^2, u_\zeta^3)) = \Span(1, u_\xi, u_\xi^2, u_\xi^3)$. Because an 
algebra homomorphism from $K \langle u_\xi \rangle = \Span(1, u_\xi, u_\xi^2, u_\xi^3)$ to~$K$ must map~$u_\xi$ to a fourth root of unity, there are four such homomorphisms, which we denote by $\alpha_\xi$, $\beta_\xi$, $\gamma_\xi$, and $\varepsilon_\xi$. By definition, they satisfy the requirements
\[\alpha_\xi(u_\xi) = \iota \qquad \beta_\xi(u_\xi) = -1 \qquad \gamma_\xi(u_\xi) = -\iota \qquad
\varepsilon_\xi(u_\xi) = 1\]
which imply that they form a basis of $K \langle u_\xi \rangle^*$. Note that $\beta_\xi$ is the restriction of the character~$\chi'_1$ considered in~\cite{KaSo3}, Par.~4.1, p.~209, while~$\varepsilon_\xi$ is the restriction of the counit~$\varepsilon_{B'}$. We have
\[\alpha_\xi(u_\xi^2) = \gamma_\xi(u_\xi^2) = -1 \qquad  \text{and} \qquad 
\beta_\xi(u_\xi^2) = \varepsilon_\xi(u_\xi^2) = 1\]

Using~$\zeta$ instead of~$\xi$, we get four algebra homomorphisms $\alpha_\zeta$, $\beta_\zeta$, $\gamma_\zeta$, and $\varepsilon_\zeta$ from $K \langle u_\zeta \rangle$ to~$K$, and since $\alpha_\xi \circ f$, $\beta_\xi \circ f$, $\gamma_\xi \circ f$, and $\varepsilon_\xi \circ f$ are also such algebra homomorphisms, these two sets must coincide.

\item
Now $\varepsilon_\xi \circ f = \varepsilon_\zeta$, as~$f$ preserves the counit.
Since $f(u_\zeta^2) = u_\xi^2$ by the first assertion, we must have $\beta_\xi \circ f = \beta_\zeta$. 
Therefore, we either have
\[\alpha_\xi \circ f = \alpha_\zeta \qquad \text{and} \qquad 
\gamma_\xi \circ f = \gamma_\zeta\]
or 
\[\alpha_\xi \circ f = \gamma_\zeta \qquad \text{and} \qquad 
\gamma_\xi \circ f = \alpha_\zeta\]

\item
In the first case, we have
\[\alpha_\xi(f(u_\zeta)) = \alpha_\zeta(u_\zeta) = \iota = \alpha_\xi(u_\xi)\]
Similarly, we get $\beta_\xi(f(u_\zeta)) = \beta_\xi(u_\xi)$, and also
$\gamma_\xi(f(u_\zeta)) = \gamma_\xi(u_\xi)$ as well as~\mbox{$\varepsilon_\xi(f(u_\zeta)) = \varepsilon_\xi(u_\xi)$.} Therefore we must have~$f(u_\zeta) = u_\xi$, as all four algebra homomorphisms in~$K \langle u_\xi \rangle ^*$ yield the same value on these two elements.

\item
In the second case, we have
\[\alpha_\xi(f(u_\zeta)) = \gamma_\zeta(u_\zeta) = -\iota = \alpha_\xi(u_\xi^3)\]
Similarly, we get $\beta_\xi(f(u_\zeta)) = \beta_\xi(u_\xi^3)$, and also
$\gamma_\xi(f(u_\zeta)) = \gamma_\xi(u_\xi^3)$ as well as \mbox{$\varepsilon_\xi(f(u_\zeta)) = \varepsilon_\xi(u_\xi^3)$.} Therefore we must have~$f(u_\zeta) = u_\xi^3$ for the same reason as in the first case.
\qed
\end{pflist}
\end{pf}

\subsection[The induced biproduct decomposition]{} \label{BiprodDec}
There is no reason why the hypothetical isomorphism~$f$ should preserve the Radford biproduct decomposition in the sense that it maps~$A \deq A_\zeta$ to \mbox{$A' \deq A_\xi$} and~$H$ to~$H$. Rather, the biproduct decomposition of~$B$ will be transported to a potentially different biproduct decomposition of~$B'$. Alternatively, we can use the inverse of~$f$ to transport the biproduct decomposition of~$B'$ to~$B$. In view of the fact that~$A' \cong B'^{\co H}$ discussed at the end of Paragraph~\ref{RadfBiprod}, we denote the factors of this decomposition by
\[\tilde{A} \deq f^{-1}(B'^{\co H}) \qquad \text{and} \qquad \tilde{H} \deq f^{-1}(H)\]
Note that this space~$\tilde{H}$ is not necessarily equal to the one denoted by the same symbol in Paragraph~\ref{GroupLikeQuot}. With this notation, we can prove the following basic fact:
\begin{prop}
If $b \in \tilde{A}$ has the property that its residue class $\overline{b} \in B/B\tilde{H}^+$ is group-like, then $f(b) \in \{c'_1, c'_2, c'_3, c'_4, d'_1, d'_2, d'_3, d'_4\}$.
\end{prop}
\begin{pf}
\begin{pflist}
\item
The isomorphism~$f$ induces a coalgebra isomorphism
\[\overline{f} \coln B/B\tilde{H}^+ \to B'/B'H{}^+\]
that maps~$\overline{b}$ to~$\overline{f(b)}$. Therefore 
$\overline{f(b)}$ is group-like in $B'/B'H{}^+$. 

\item
By hypothesis, we have $f(b) \in B'^{\co H}$. As just recalled,~$A' \cong B'^{\co H}$, so there exists an element~$a' \in A'$ with~$f(b) = a' \star \HH$, which obviously implies that~$\overline{f(b)} = \overline{a' \star \HH}$. Now Proposition~\ref{RadfBiprod} implies that~$a'$ is group-like in~$A'$, so that~$a'= \omega'_i$ or~$a' = \eta'_i$ for some~$i=1,2,3,4$ (cf.~\cite{Sweed}, Prop.~3.2.1, p.~54). But this means that~$f(b) = c'_i$ or~$f(b) = d'_i$ for some~$i=1,2,3,4$, as asserted.~\qed
\end{pflist}
\end{pf}

We will in particular be interested in the following special case:
\begin{corollary}
Suppose that $b \in \tilde{A}$ is not contained in $\Span(1, u_\zeta, u_\zeta^2, u_\zeta^3)$. If the residue class $\overline{b} \in B/B\tilde{H}^+$ is group-like, then $f(b) \in \{d'_1, d'_2, d'_3, d'_4\}$ and the order of~$b$ is equal to the order of~$v_\xi = d'_1$.
\end{corollary}
\begin{pf}
As already pointed out in the proof of Proposition~\ref{IsomFirstZeta}, it follows from~\cite{KaSo3}, Par.~2.1, p.~203 that 
$\Span(1, u_\xi, u_\xi^2, u_\xi^3) = \Span(c'_1, c'_2, c'_3, c'_4)$. Since 
\[f(\Span(1, u_\zeta, u_\zeta^2, u_\zeta^3)) = \Span(1, u_\xi, u_\xi^2, u_\xi^3) \]
by Proposition~\ref{IsomBiprodZetaOne}, we then must have
\[f(b) \in \{d'_1, d'_2, d'_3, d'_4\} = \{v_\xi, u_\xi v_\xi, u_\xi^2 v_\xi, u_\xi^3 v_\xi\}\]
Since $u_\xi$ has order~$4$, all elements in the set $\{v_\xi, u_\xi v_\xi, u_\xi^2 v_\xi, u_\xi^3 v_\xi \}$ have the same order, namely the order of~$v_\xi$, which is~$4$ if~$\xi$ is primitive and~$8$ if~$\xi$ is not primitive, as recalled in Paragraph~\ref{IsomFirstZetaOne}. This shows that~$f(b)$, and therefore also~$b$, has this order.
\qed
\end{pf}

\subsection[The image of~$s$]{} \label{s}
We can restrict the possible images of the generators further:
\begin{prop} 
$f(s_\zeta) = s_\xi$
\end{prop}
\begin{pf}
By~\cite{KaSo3}, Thm.~6.3, p.~218, $B$ contains a unique four-dimensional normal Hopf subalgebra, which is spanned by the group-like elements~$\B$, $u_\zeta^2$, $s_\zeta$, and~$u_\zeta^2 s_\zeta$. The isomorphism~$f$ must map this algebra to the corresponding normal Hopf subalgebra of~$B'$. Since it also maps group-like elements to group-like elements, it follows from Proposition~\ref{IsomBiprodZetaOne} that $f(s_\zeta) = s_\xi$ or~$f(s_\zeta) = u_\xi^2 s_\xi$. We want to rule out the second possibility.

So suppose that~$f(s_\zeta) = u_\xi^2 s_\xi$. Then we also have 
$f(u_\zeta^2 s_\zeta) = u_\xi^2 f(s_\zeta) = s_\xi$ by Proposition~\ref{IsomBiprodZetaOne}.
With the notation~$\tilde{A}$ and~$\tilde{H}$ introduced in Paragraph~\ref{BiprodDec}, we get 
\[\tilde{H} = K \langle f^{-1}(r_\xi), u_\zeta^2 s_\zeta \rangle\]
We then have $u_\zeta (u_\zeta^2 s_\zeta - \B) \in B\tilde{H}{}^+$ and consequently
$\overline{u_\zeta^3 s_\zeta} = \overline{u}_\zeta$ in the quotient coalgebra~$B/B\tilde{H}{}^+$. The coproduct formula 
\begin{align*}
\db(u_\zeta) &= \frac{\K}{2} (u_\zeta \ot u_\zeta + u_\zeta \ot u_\zeta^3 +  u_\zeta^3 s_\zeta \ot u_\zeta - u_\zeta^3 s_\zeta \ot u_\zeta^3)
\end{align*}
recalled in Paragraph~\ref{RadfBiprodCase} therefore yields 
\begin{align*}
\Delta(\overline{u}_\zeta) &= \frac{\K}{2} (\overline{u}_\zeta \ot \overline{u}_\zeta + \overline{u}_\zeta \ot \overline{u_\zeta^3} 
+  \overline{u}_\zeta \ot \overline{u}_\zeta - \overline{u}_\zeta \ot \overline{u_\zeta^3})
= \overline{u}_\zeta \ot \overline{u}_\zeta
\end{align*}
so $\overline{u}_\zeta$ is group-like in~$B/B\tilde{H}{}^+$. From Proposition~\ref{IsomBiprodZetaOne}, we know that $f(u_\zeta) = u_\xi$ or $f(u_\zeta) = u_\xi^3$, which implies that 
$u_\zeta \in \tilde{A}$. But now Proposition~\ref{BiprodDec} yields 
\[f(u_\zeta) \in \{c'_1, c'_2, c'_3, c'_4, d'_1, d'_2, d'_3, d'_4\}\]
so that we have reached a contradiction.
\qed
\end{pf}

As a consequence, we can also restrict the possible values for~$f(r_\zeta)$:
\begin{corollary}
We have
$f(r_\zeta) = r_\xi$, $f(r_\zeta) = r_\xi s_\xi$, $f(r_\zeta) = u_\xi^2 r_\xi$, 
or $f(r_\zeta) = u_\xi^2 r_\xi s_\xi$.
\end{corollary}
\begin{pf}
Since~$r_\zeta$ is group-like, $f(r_\zeta)$ must also be group-like. We have already mentioned in Paragraph~\ref{GroupLikeQuot} that~$G(B_\xi) = \langle u_\xi^2, r_\xi, s_\xi \rangle$, which contains eight elements. Since $f(\B) = 1_{B'}$, $f(u_\zeta^2) = u_\xi^2$, $f(s_\zeta) = s_\xi$, and $f(u_\zeta^2 s_\zeta) = u_\xi^2 s_\xi$, only the four remaining elements are possible as values of~$f(r_\zeta)$.
\qed
\end{pf}

\subsection[Subgroups of~$G(B^*)$]{} \label{Subgr}
In Paragraph~\ref{PropRadfBiprod}, we have considered the generators~$\chi_1$,~$\chi_2$, and~$\chi_3$ of~$G(B^*) = G(B_\zeta^*)$. We denote the analogous characters of~$G(B'^*) = G(B_\xi^*)$ by~$\chi'_1$,~$\chi'_2$, and~$\chi'_3$, which therefore need to be distinguished from the characters considered in~\cite{KaSo3}, Par.~4.1, p.~209, which were denoted by the same symbols and were briefly mentioned in Paragraph~\ref{IsomBiprodZetaOne}.

The group $\Gamma \deq f^*(\langle \chi'_2, \chi'_3 \rangle)$ is a subgroup of order~$4$ in~$G(B^*_\zeta)$. The same argument as in~\cite{KaSo3P}, Par.~6.3, p.~36 shows that there are seven such subgroups, namely \mbox{$\Gamma_1 \deq \langle \chi_1, \chi_2 \rangle = \{\eb, \chi_1, \chi_2, \chi_1 \chi_2\}$} and
\begin{align*}
\Gamma_2 &\deq \{\eb, \chi_1, \chi_3, \chi_1 \chi_3\} & 
\Gamma_3 &\deq \{\eb, \chi_1, \chi_2 \chi_3, \chi_1 \chi_2 \chi_3\} \\
\Gamma_4 &\deq \{\eb, \chi_2, \chi_3, \chi_2 \chi_3\} & \Gamma_5 &\deq \{\eb, \chi_2, \chi_1 \chi_3, \chi_1 \chi_2 \chi_3\} \\
\Gamma_6 &\deq  \{\eb, \chi_3, \chi_1 \chi_2, \chi_1 \chi_2 \chi_3\}&
\Gamma_7 &\deq \{\eb, \chi_2 \chi_3, \chi_1 \chi_3, \chi_1 \chi_2\}
\end{align*}
By~\cite{KaSo3}, Prop.~3.1, p.~207, the character~$\chi'_1$ is the unique nontrivial central group-like element in~$G(B_\xi^*)$. Since $\langle \chi'_2, \chi'_3 \rangle$ does not contain~$\chi'_1$, the image $\Gamma = f^*(\langle \chi'_2, \chi'_3 \rangle)$ cannot contain~$\chi_1$ and therefore cannot be equal to~$\Gamma_1$, $\Gamma_2$, or~$\Gamma_3$. As we will show now, this image can also not be equal to one of the remaining groups if either~$\zeta$ or~$\xi$ is primitive and the other is not, thereby completing the proof of Theorem~\ref{IsomBiprodZetaOne}.

As we stated already in the proof of Corollary~\ref{RadfBiprod}, the map
\[A \to B^{\co H},~a \mapsto a \star \HH\]
is bijective, which implies that the algebra~$B^{\co H}$ is generated by~$u_\zeta$ and~$v_\zeta$. It therefore follows from Proposition~\ref{PropRadfBiprod} that
\begin{align*}
B^{\co H} &= \{b \in B \mid \psi_2(b) = b \; \text{and} \; \psi_3(b) = b\} \\
&= \{b \in B \mid \psi_\chi(b) = b \; \text{for all }  \chi \in \langle \chi_2, \chi_3 \rangle\}
\end{align*}
A similar statement holds for~$B'=B_\xi$. We therefore have
\[f^{-1}(B_\xi^{\co H}) = 
\{b \in B_\zeta \mid \psi_\chi(b) = b \; \text{for all } \chi \in \Gamma \}\]
Using this observation, it is not complicated to rule out two more of the remaining groups:
\begin{parlist}
\item
If $\Gamma = \Gamma_6 = \langle \chi_3, \chi_1 \chi_2 \rangle$, we have
\[f^{-1}(B_\xi^{\co H}) = 
\{b \in B_\zeta \mid \psi_3(b) = b \; \text{and} \; \psi_1(\psi_2(b)) = b \}\]
Then Proposition~\ref{PropRadfBiprod} implies that $f^{-1}(B_\xi^{\co H})$ contains the elements
$u_\zeta$ and~$v_\zeta r_\zeta$. But according to the relations recalled in Paragraph~\ref{RadfBiprodCase}, these two elements do not commute. On the other hand, $f^{-1}(B_\xi^{\co H})$ is clearly isomorphic to the commutative algebra~$B_\xi^{\co H}$. Therefore, this case cannot occur.

\item 
If $\Gamma = \Gamma_7 = \langle \chi_1 \chi_2, \chi_1 \chi_3 \rangle$, we have
\[f^{-1}(B_\xi^{\co H}) = 
\{b \in B_\zeta \mid \psi_1(\psi_2(b)) = b \; \text{and} \; \psi_1(\psi_3(b)) = b \}\]
As in the previous case, Proposition~\ref{PropRadfBiprod} implies that $f^{-1}(B_\xi^{\co H})$ contains the elements $u_\zeta$ and~$v_\zeta r_\zeta s_\zeta$. Again according to the relations recalled in Paragraph~\ref{RadfBiprodCase}, these two elements do not commute, so that this case can also not occur.
\end{parlist}

\subsection[The group~$\Gamma_4$]{} \label{Gamma4}
If $\Gamma = \Gamma_4 = \langle \chi_2, \chi_3 \rangle$, we have
\[f^{-1}(B_\xi^{\co H}) = 
\{b \in B_\zeta \mid \psi_2(b) = b \; \text{and} \; \psi_3(b) = b \}\]
As discussed in Paragraph~\ref{Subgr}, this means that $f^{-1}(B_\xi^{\co H}) = B_\zeta^{\co H} = \langle u_\zeta, v_\zeta \rangle$, so that in particular
$f(v_\zeta) \in B_\xi^{\co H} = \langle u_\xi, v_\xi \rangle$. From Corollary~\ref{s}, we know that $f(r_\zeta) = r_\xi$, $f(r_\zeta) = r_\xi s_\xi$, $f(r_\zeta) = u_\xi^2 r_\xi$, or 
$f(r_\zeta) = u_\xi^2 r_\xi s_\xi$. We look at the cases separately:
\begin{parlist}
\item 
First, we assume that $f(r_\zeta) = r_\xi$ or $f(r_\zeta) = r_\xi s_\xi$. Then we get from Proposition~\ref{s} that $f(H) = H$, and consequently $f^{-1}(H) = H$.
The residue class of~$v_\zeta$ in $B/BH{}^+$ is group-like, as we noted at the beginning of Paragraph~\ref{GroupLikeQuot}. Now Corollary~\ref{BiprodDec} implies that
$v_\zeta$ has the same order as~$v_\xi$, which is~$4$ if~$\xi$~is primitive and~$8$ if~$\xi$~is not primitive, as mentioned in the proof of that corollary. This implies that either both~$\xi$ and~$\zeta$ are primitive or both are not primitive. 

\item
Let us now consider the case $f(r_\zeta) = u_\xi^2 r_\xi$ or $f(r_\zeta) = u_\xi^2 r_\xi s_\xi$. Using the Hopf subalgebra
$\tilde{H} \deq \Span(\langle  u_\zeta^2 r_\zeta, s_\zeta \rangle)$ from Paragraph~\ref{GroupLikeQuot}, we then have 
$f(\tilde{H}) = H$, as required in Paragraph~\ref{BiprodDec}, since
$f(u_\zeta^2 r_\zeta) = u_\xi^2 (u_\xi^2 r_\xi) = r_\xi$ or \mbox{$f(u_\zeta^2 r_\zeta) = u_\xi^2 (u_\xi^2 r_\xi s_\xi) = r_\xi s_\xi$}, as well as $f(s_\zeta) = s_\xi$ by Proposition~\ref{s}.
From Proposition~\ref{GroupLikeQuot}, we know that the residue class of~$\frac{1+\iota}{2} v_\zeta 
+ \frac{1-\iota}{2} u_\zeta^2 v_\zeta$ is group-like in $B/B\tilde{H}{}^+$. Now Corollary~\ref{BiprodDec} implies that 
$\frac{1+\iota}{2} v_\zeta + \frac{1-\iota}{2} u_\zeta^2 v_\zeta$ has the same order as~$v_\xi$. But we have
\begin{align*}
\left(\frac{1+\iota}{2} v_\zeta + \frac{1-\iota}{2} u_\zeta^2 v_\zeta \right)^2 &= 
\frac{\iota}{2} v_\zeta^2 + v_\zeta (u_\zeta^2 v_\zeta) 
- \frac{\iota}{2} (u_\zeta^2 v_\zeta)^2 = u_\zeta^2 v_\zeta^2 
\end{align*}
and therefore
\begin{align*}
\left(\frac{1+\iota}{2} v_\zeta + \frac{1-\iota}{2} u_\zeta^2 v_\zeta \right)^4 &=  v_\zeta^4
\end{align*}
This shows that $\frac{1+\iota}{2} v_\zeta + \frac{1-\iota}{2} u_\zeta^2 v_\zeta$ has not only the same order as~$v_\xi$, but also the same order as~$v_\zeta$, so that~$v_\xi$ and~$v_\zeta$ have the same order. As in the previous case, this implies that either both~$\xi$ and~$\zeta$ are primitive or both are not primitive. 
\end{parlist}

\subsection[The group~$\Gamma_5$]{} \label{Gamma5}
If $\Gamma = \Gamma_5 = \langle \chi_2, \chi_1 \chi_3 \rangle$, we have
\[f^{-1}(B_\xi^{\co H}) = 
\{b \in B_\zeta \mid \psi_2(b) = b \; \text{and} \; \psi_1(\psi_3(b)) = b \} \]
As before, we now apply Proposition~\ref{PropRadfBiprod}, which in the present case yields that~$f^{-1}(B_\xi^{\co H})$ contains the elements~$u_\zeta$ and~$v_\zeta s_\zeta$. From Corollary~\ref{s}, we know that $f(r_\zeta) = r_\xi$, \mbox{$f(r_\zeta) = r_\xi s_\xi$}, $f(r_\zeta) = u_\xi^2 r_\xi$, or $f(r_\zeta) = u_\xi^2 r_\xi s_\xi$, and again we look at the cases separately:

\begin{parlist}
\item
As in the treatment of~$\Gamma_4$ in Paragraph~\ref{Gamma4}, we assume first that $f(r_\zeta) = r_\xi$ or $f(r_\zeta) = r_\xi s_\xi$. Then we have \mbox{$f(H) = H$} by Proposition~\ref{s}, and consequently $f^{-1}(H) = H$. The residue class of $v_\zeta s_\zeta$ in $B/BH{}^+$ is group-like, as we noted at the beginning of Paragraph~\ref{GroupLikeQuot}. Now Corollary~\ref{BiprodDec} implies that
$v_\zeta s_\zeta$ has the same order as~$v_\xi$, which is~$4$ if~$\xi$ is primitive and~$8$ if~$\xi$ is not primitive, as just recalled again in Paragraph~\ref{Gamma4}. 

Using the relations stated in Paragraph~\ref{RadfBiprodCase}, we have
\[(v_\zeta s_\zeta)^2 = v_\zeta (s_\zeta v_\zeta) s_\zeta = u_\zeta^2 v_\zeta^2 \]
and so $(v_\zeta s_\zeta)^4 = v_\zeta^4$. Thus if one of~$v_\xi$ and~$v_\zeta$ has order~$4$, then so does the other. This implies that either both~$\xi$ and~$\zeta$ are primitive or both are not primitive. 

\item
Let us now consider the case $f(r_\zeta) = u_\xi^2 r_\xi$ or $f(r_\zeta) = u_\xi^2 r_\xi s_\xi$. The same computation as in the corresponding case in Paragraph~\ref{Gamma4} then shows that the Hopf subalgebra $\tilde{H} \deq \Span(\langle u_\zeta^2 r_\zeta, s_\zeta \rangle)$ from Paragraph~\ref{GroupLikeQuot} satisfies $f(\tilde{H}) = H$, as required in Paragraph~\ref{BiprodDec}. From the above, we have
\[\frac{1+\iota}{2} v_\zeta s_\zeta + \frac{1-\iota}{2} u_\zeta^2 v_\zeta s_\zeta
 \in f^{-1}(B_\xi^{\co H})\]
and as before Proposition~\ref{GroupLikeQuot} yields that the residue class of this element is group-like in $B/B\tilde{H}{}^+$. So again Corollary~\ref{BiprodDec} implies that this element has the same order as~$v_\xi$. Since we have already computed that $(v_\zeta s_\zeta)^2 = u_\zeta^2 v_\zeta^2$, we have
\begin{align*}
\left(\frac{1+\iota}{2} v_\zeta s_\zeta + \frac{1-\iota}{2} u_\zeta^2 v_\zeta s_\zeta \right)^2 &= 
\frac{\iota}{2} (v_\zeta s_\zeta)^2 + (v_\zeta s_\zeta) (u_\zeta^2 v_\zeta s_\zeta) 
- \frac{\iota}{2} (u_\zeta^2 v_\zeta s_\zeta)^2 \\
&= 
(v_\zeta s_\zeta) (u_\zeta^2 v_\zeta s_\zeta) = v_\zeta^2
\end{align*}
and therefore
\begin{align*}
\left(\frac{1+\iota}{2} v_\zeta s_\zeta + \frac{1-\iota}{2} u_\zeta^2 v_\zeta s_\zeta \right)^4 &=  v_\zeta^4
\end{align*}
This shows that $\frac{1+\iota}{2} v_\zeta s_\zeta + \frac{1-\iota}{2} u_\zeta^2 v_\zeta s_\zeta$ has not only the same order as~$v_\xi$, but also the same order as~$v_\zeta$, so that~$v_\xi$ and~$v_\zeta$ have the same order. As in the previous case, this implies that either both~$\xi$ and~$\zeta$ are primitive or both are not primitive. And as we have now considered all possibilities for the subgroup~$\Gamma$, this finishes the proof of Theorem~\ref{IsomBiprodZetaOne}.
\end{parlist}

\end{document}